\documentclass[11pt]{article}

\usepackage{setspace,graphicx,amssymb,amsmath,mathtools,latexsym,amsfonts,amscd,amsthm,multirow,mathdots,caption,array,dsfont,wrapfig}
\usepackage{fancyhdr,tabularx,cite,mathrsfs,enumitem,graphicx}
\usepackage{stmaryrd}
\usepackage{rotating}
\usepackage{authblk}
\usepackage{hyperref}
\usepackage{comment,enumitem}
\usepackage{geometry}
\usepackage{arydshln}
\usepackage{longtable}
\usepackage{hyperref}
\usepackage{booktabs}
\usepackage{braket}
\usepackage{subcaption} 
\usepackage{lineno}

\usepackage[font=small,labelfont=bf]{caption}
\usepackage{tikz}
\usetikzlibrary{calc,trees,positioning,arrows,chains,shapes.geometric,%
    decorations.pathreplacing,decorations.pathmorphing,shapes,%
    matrix,shapes.symbols,decorations.markings}
\usepackage{float}

\newcommand{\bea}{\begin{eqnarray}} 
\newcommand{\eea}{\end{eqnarray}} 
\newcommand{\bee}{\begin{eqnarray*}} 
\newcommand{\eee}{\end{eqnarray*}} 
\newcommand{\al}{\begin{align*}} 
\newcommand{\eal}{\end{align*}} 
\newcommand{\be}{\begin{equation}} 
\newcommand{\ee}{\end{equation}} 
\newcommand{\eq}[1]{(\ref{#1})} 
\newcommand{\bem}{\begin{pmatrix}} 
\newcommand{\eem}{\end{pmatrix}} 

\newcommand{\GS}[1]{{\color{green}{GS: #1}}}
\newcommand{\FF}[1]{{\color{blue}{FF: #1}}}
 
\def\c{\gamma}

\def\t{\tau}

\newcolumntype{R}{ >{$}r <{$}}
\newcolumntype{C}{ >{$}c <{$}}
\newcolumntype{L}{ >{$}l <{$}}
\newcolumntype{F}{>{\centering\arraybackslash}m{1.5cm}}

\newcommand{\RR}{{\mathbb R}}
\newcommand{\CC}{{\mathbb C}}
\newcommand{\ZZ}{{\mathbb Z}}
\newcommand{\QQ}{{\mathbb Q}}
\newcommand{\HH}{{\mathbb H}}

\newcommand{\SL}{\operatorname{\textsl{SL}}}      

\newcommand{\g}{\gamma}	

\newcommand{\xmod}{{\rm \;mod\;}}

\newcommand{\CS}{\rm CS}

\newtheorem{thm}{Theorem}[section]

\newtheorem{lem}[thm]{Lemma}
\newtheorem{prop}[thm]{Proposition}
\newtheorem{conj}[thm]{Conjecture}

\theoremstyle{definition}
\newtheorem{defn}[thm]{Definition}

\theoremstyle{remark}

\numberwithin{equation}{section}
\newtheorem*{eg}{Example}

\title{
\vspace{-35pt}
    \textsc{\Large{
Three-Manifold Quantum Invariants and Mock Theta Functions
            }  }
}

\author[1,2]{\small{Miranda C. N. Cheng}\thanks{On leave from CNRS, France.}}
\author[3,4]{\small{Francesca Ferrari}}
\author[2]{\small{Gabriele Sgroi}\thanks{g.sgroi@uva.nl}}

\date{}

\affil[1]{Korteweg-de Vries Institute for Mathematics, Amsterdam, the Netherlands}
\affil[2]{Institute of Physics, University of Amsterdam, Amsterdam, the Netherlands}
\affil[3]{International School for Advanced Studies (SISSA), Trieste, Italy} 
\affil[4]{INFN, Sezione di Trieste, Trieste, Italy}

\begin{document}

\setstretch{1.4}

   \maketitle
   \abstract
   {Mock modular forms have found applications in numerous branches of mathematical sciences since they were first introduced by Ramanujan nearly a century ago. In this proceeding we highlight a new area where mock modular forms start to play an important role, namely the study of three-manifold invariants. 
   For a certain class of Seifert three-manifolds,
we describe a conjecture on the mock modular properties of a recently proposed quantum invariant. 
  As an illustration, we include concrete computations for a specific  three-manifold,  the Brieskorn sphere $\Sigma(2,3,7).$  
  This note is partially based on the talk by the first author in the conference ``Srinivasa Ramanujan: in celebration of the centenary of his election as FRS'' held at the Royal Society in 2018. 
   }

\newpage 
\tableofcontents

\section{Introduction}\label{intro}

Mock modular forms have been a source of fascination since Ramanujan first introduced them about a century ago \cite{Ram88},\cite{Ram00}. 
In 2002, Zweger's doctoral thesis provided the crucial starting point for a structural theory of mock modular forms. 
These functions have found applications in various branches of mathematical sciences, including combinatorics, moonshine, conformal field theory, string theory and more. 
We refer to \cite{zagier_mock}, \cite{Ono_unearthing}, \cite{Folsom_what}, \cite{Duk_almost},\cite{UM}, \cite{MUM}, \cite{book} and other contributions to this volume for a partial overview of these developments. 
The purpose of this article is to highlight  the appearance of mock modular forms  in a different context: the topology of three-dimensional manifolds. 
This appearance was first anticipated in \cite{3d} and later further developed via concrete computations in \cite{GM} and \cite{future}.


To see how mock modular forms appear in the study of three-manifolds, we first introduce a set of topological invariants, noted by $\widehat Z_a$, defined in \cite{BPS} for {\it weakly negative plumbed manifolds}. 
Roughly speaking, these are three-manifolds obtained through surgeries along links that are in turn determined by weighted graphs (cf. Figure \ref{fig:plumbing}), which moreover satisfy a certain negativity condition\cite {3d}.

More precisely, the data we need is a plumbing graph, which is a 
 a weighted simple graph $(V,E,\lambda)$ specified by 
the set $V$ of vertices, the set $E$ of edges, and an integral weight function $\lambda:V\to \ZZ$. 
Equivalently, the data can be captured by an adjacency matrix $M$, which is a square matrix of size $|V|$ with entries given by  $M_{vv'}= \lambda(v)$ if $v=v'$, $1$ if $(v,v')\in E$ and 0 otherwise. The data determines a three-manifold\footnote{Different weighted graph related by the so-called Kirby moves can lead to the same topological three-manifold $M_3$. See \cite {GM} for a proof the Kirby-invariance of the quantum invariants $\widehat Z_a(M_3)$ defined in Definition \ref{def:hatZ}.}. We say that $M_3$ is a weakly negative plumbed manifold if $M^{-1}$ is negative-definite when restricted to the subspace generated by all vertices with degree larger than 2. 

\begin{defn}\label{def:hatZ}
For $M_3$ a weakly negative plumbed three-manifold and using the above notation, we define 
 the quantum invariants $\widehat{Z}_a(M_3;\tau)$ via the following principal value $|V|$-dimensional  integral: 
 \be
\widehat{Z}_a(M_3;\tau):= (-1)^{\pi}q^{\frac{3\sigma -\sum_{v\in V} \lambda(v)}{4}} \,\sum_{{\bf n}\in 2M\ZZ^{|V|}+{\bf a}} \text{vp} \prod_{v\in V} \oint_{|w_v|=1}\frac{dw_v}{2\pi i w_v}\biggl(w_v-\frac{1}{w_v}\biggr)^{2-\text{deg}(v)}q^{-\frac{{\bf n}^TM^{-1}{\bf n}}{4}}\  e^{2\pi i{\bf z}^T {\bf n}}
\ee
 where we write $q:=e^{2\pi i \tau}$ and $w_v:=e^{2\pi i z_v}$ as usual, and use the bold-faced letters to denote elements in $\ZZ^{|V|}$. When $M^{-1}$ is moreover negative definite, the above can be rewritten as
\be
\widehat{Z}_a(M_3;\tau):= (-1)^{\pi}q^{\frac{3\sigma -\sum_{v\in V} \lambda(v)}{4}} \, \text{vp} \prod_{v\in V} \oint_{|w_v|=1}\frac{dw_v}{2\pi i w_v}\biggl(w_v-\frac{1}{w_v}\biggr)^{2-\text{deg}(v)}\Theta^{-M}_a(\tau, {\bf z})
\ee
In the above, $\pi$ denotes the number of positive  eigenvalues, and $\sigma$ is the signature of $M^{-1}$.
The label $a$ of the quantum invariants $\widehat{Z}_a(M_3)$ can be identified with elements of the set ${\rm Spin}^c(Y)\cong  \pi_0{\cal M}_{\rm ab}(M_3)\cong (2\ZZ^{|V|}+\delta)/(2M\ZZ^{|V|})$, where
$\delta \in \ZZ^{|V|}/2\ZZ^{|V|}$ is defined by
{$\delta_v={\rm deg}(v)$ mod 2}, and ${\cal M}_{\rm ab}(M_3)$ denotes the moduli space of Abelian flat connections. Denote by ${\bf a}$ the corresponding element of $(2\ZZ^{|V|}+\delta)/(2M\ZZ^{|V|})$,  the theta function reads
\be\label{def:theta}
\Theta^{-M}_a(\tau, {\bf z}) = \sum_{{\bf n}\in 2M\ZZ^{|V|}+{\bf a}} q^{-\frac{{\bf n}^TM^{-1}{\bf n}}{4}}\  e^{2\pi i{\bf z}^T {\bf n}}.
\ee 
\end{defn}

\vspace{2pt}
\newsavebox\mybox
	\begin{lrbox}{\mybox}
		\tikzstyle{point}=[circle,draw=black!100,fill=black!100,thick,
		inner sep=0pt,minimum size=2mm]
		\begin{tikzpicture}[scale=0.8]
	\node[point] (a1) at (0,1)  [label=left:{\small $v_1$}]{};
	\node[point] (a2) at (0,-1)[label=left:{\small $v_2$}] {};
	\node[point] (a3) at (1,0)[label=left:{\small$v_3$}] {};
	\node[point] (a4) at (3,0)[label=right:{\small $v_4$}] {};
	\node[point] (a5) at (4,1)[label=right:{\small $v_5$}] {};
	\node[point] (a6) at (4,-1)[label=right:{\small $v_6$}] {};
	\draw (a1)--(a3)--(a2); \draw (a3)--(a4); \draw (a5)--(a4)--(a6); 
		\end{tikzpicture}
	\end{lrbox}
	\begin{figure}
		\centering
	\begin{tikzpicture}[every text node part/.style={align=center}]
	\node (A) at (0,0)[label={[label distance=0.5cm]above:{Plumbing graph $\Gamma$,\\with weights $\lambda:V\to \ZZ$}}]{
		\usebox\mybox
	};
\node (B) at (7,0)[label={[label distance=0.5cm]above:{Framed link $L$,\\ with framing coefficients $\lambda(v_i)$}}]{
	\includegraphics[scale=0.24]{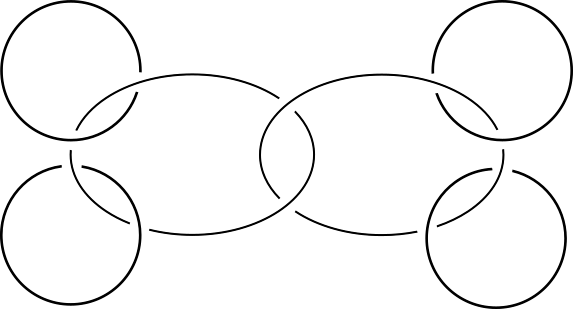}
};
\node (C) at (0,-5)[label=above:{Adjacency matrix $M$}]{
	$
	\begin{pmatrix}
	\lambda(v_1) & 0 & 1 &0 & 0 &0 \\
	0 & \lambda(v_2) &1 & 0 & 0 &0 \\
	1 & 1 & \lambda(v_3) & 1 & 0 & 0 \\
	0 & 0 & 1 & \lambda(v_4) & 1 & 1 \\
	0 & 0 & 0 & 1 & \lambda(v_5) & 0 \\
	0 & 0 & 0 & 1 & 0 & \lambda(v_6)
	\end{pmatrix}$
	}; 
\node(D) at (7,-4) {Plumbed $M_{3}$ \\ obtained from surgery along $L$};
\draw[<->,line width=.2mm, shorten >=6pt] ($(A)+(2.5,0)$)--($(B)+(-2.5,0)$);
\draw[->,line width=.2mm, shorten >=6pt] ($(B)+(0,-1.5)$)--($(D)+(0,.5)$);\draw[<->,line width=.2mm, shorten >=6pt] ($(C)+(0,2.5)$)--($(A)+(0,-1)$);
	\end{tikzpicture}
	\caption{\label{fig:plumbing}Weighted graphs, adjacency matrices, links and plumbed three-manifolds. }
\end{figure}
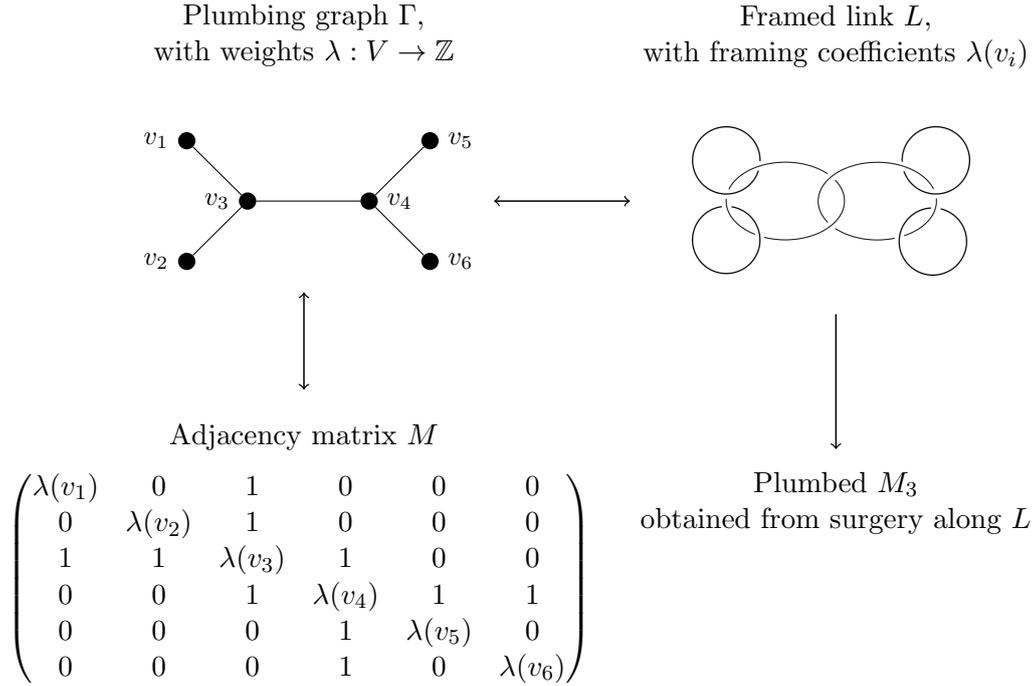

A well-known  topological invariant for three-manifolds is the 
Witten-Reshetikhin-Turaev 
(WRT) invariant, defined for all three-manifolds. Physically speaking, it is (up to certain well-understood prefactors) the  partition function of Chern-Simons theory on the three-manifold $M_3$ which we denote by $Z_{\rm CS}(M_3)$. For a given three-manifold $M_3$ (and a simple Lie group $G$ which we will take to be $G=SU(2)$ for the sake of concreteness), we obtain a function $Z_{\rm CS}(M_3): \ZZ \to \CC$ defined on all integers, namely the (shifted) Chern-Simons levels. Analogous to knot theory, it would be desirable to have a $q$-series version of the invariants defined on a continuous domain, such as the upper-half plane. This would then be the first step towards a categorification of three-manifold invariants, analogous to the categorification programme of knot invariants. 
It was shown for weakly negative plumbed manifolds that the quantum invariants $\widehat Z_a$ provide exactly such a $q$-series generalisation of the WRT invariants. 
To be more precise, 
for weakly negative plumbed manifolds one has \cite{BPS}  
\be
\label{eq:wrt}
(i\sqrt{2k}) Z_{\CS}(M_3;k) = \sum_{a, b} X_{ab} e^{2\pi i {{\CS}(a)} k}\lim_{\tau\to \frac{1}{k}} \widehat{Z}_b(M_3;\tau).
\ee
In the above equation, the sum over $a$ is over the set $((2\ZZ^{|V|}+\delta)/2M\ZZ^{|V|}))/\ZZ_2$,  which can be identified with the space of gauge-inequivalent $SU(2)$ Abelian flat connections on $M_3$, and 
 ${\CS}(a)$ denotes the corresponding Chern-Simons invariant when we regard $a$ as a  label for Abelian $SU(2)$ flat connections. The 
 sum over $b$ is over the set $(\ZZ^{|V|}/M\ZZ^{|V|})/\ZZ_2$, and the matrix $X$ has as elements 
\begin{equation}\label{def:Xmatrix}
    X_{ab}=\frac{
    \sum\limits_{
     (a',b')\in\{\mathbb{Z}_2\times \mathbb{Z}_2~ \text{orbit of }(a,b)\} }
   e^{2\pi i (a', M^{-1}b')}}
   {2\sqrt{|\text{Det}M|} }.
\end{equation}
To summarise, two steps need to be taken in order to retrieve $Z_{\CS}$ from $\widehat Z_a$. First, $\widehat Z_a$ has an extra label $a$ indexing the $SU(2)$ Abelian flat connections while $Z_{\CS}$ does not, and this label needs therefore to be summed over. Second, a so-called radial limit $\tau \to \frac{1}{k}$ taking $\tau\in \HH$ to the boundary $\QQ\cup \{i\infty\}$ of the upper-half plane needs to be taken in order to relate the continuous variable $\tau$ and the (shifted)  Chern-Simons level.

The modular-like properties of the quantum invariants $\widehat Z_a$ is a rich subject that has been in development since \cite{3d}. So far it develops in parallel to the study of modular-like properties of knot invariants (see for instance \cite{qmf}, \cite{GZ}, \cite{DG}, \cite{HL} for a sample of work in this direction), although it is expected that the two topics are related  both in their physical and mathematical contexts. 

For concreteness and in order to make direct contact with Ramanujan's mock theta functions, here we 
 restrict our attention to the simplest non-trivial plumbing graphs: the so-called three-star weighted graphs.  These are, as the name suggests, weighted simple graphs with one vertex of degree three, three vertices  of degree one, while the rest of the vertices (if any) have degree two. See Figure \ref{fig:3leg}. We will denote the unique vertex with degree three by $v_0$. 
Such graphs are either weakly negative or not, depending on the sign of $(M^{-1})_{00}$. 
When $(M^{-1})_{00}<0$, Definition \ref{def:hatZ} is readily applicable and it is not hard to show that the quantum invariants $\widehat Z_a$  are always holomorphic functions on the upper-half plane with well-defined $q$-expansions and moreover have integral coefficients. In fact, a lot more is true: up to a possible addition of a polynomial, the quantum invariants $\widehat Z_a$ are linear combinations of false theta functions multiplied by a rational $q$-power (cf. \S\ref{subsec:false3star}  and \cite{plumbing}).

A puzzle immediately arises given the simple result for weakly negative three-star graphs: what happens when one flips the orientation of the three-manifold? While  this might sound like  an innocuous  operation, it can in fact have rather dramatic consequences due to the  {\it pseudo-chiral symmetry} (or CP symmetry in physical terms) 
\be\label{eqn:rel_CS_orientation}
Z_{\CS}(M_3;k)=Z_{\CS}(-M_3;-k)
\ee
of Chern-Simons theory. 
From the relation \eqref{eq:wrt} between the quantum invariants $\widehat Z_a(M_3;\tau)$ and $Z_{\CS}(M_3;k)$, and in particular the relation ``$\tau \to \frac{1}{k}$'' between the two variables, one is led to the guess 
\be\label{eqn:reflect}
\widehat Z_a(-M_3;\tau) {\text{ ``=''}}\, \widehat Z_a(M_3;-\tau). 
\ee
There are a few immediate problems with this guess. 
Recall that  for  a weakly negative plumbed manifold $M_3$, Definition \ref{def:hatZ} defines a function $\widehat Z_a(M_3;\tau)$ 
 on the upper-half plane $\HH$, which is not preserved by the action $\tau \mapsto -\tau$. As a result it is not clear what the right-hand side of the equation \eq{eqn:reflect} even means. More concretely, it is clear from \eq{def:hatZ} that for plumbed manifolds one has $\tau \mapsto -\tau \Leftrightarrow M \mapsto -M$, which flips the sign of the adjacency matrix and hence flips the signature of the lattice for which the theta function $\Theta^{-M}_a$ should be defined, and as a result does not render a function on $\HH$ when one tries to  literally apply Definition \ref{def:hatZ}. 
 

This is when the question starts to become interesting from the perspective of mock modular forms. To be concrete, we let $M_3$ be a weakly negative three-star plumbed three-manifold.  
As mentioned before, for such cases $\widehat Z_a(M_3;\tau)$ are basically false theta functions, which are known to furnish (rather simple) examples of the so-called {\it quantum modular forms}, as will be explained in \S\ref{sec:false_mock}. The quantum modular properties of the quantum invariants $\widehat Z_a$ are essentially what makes their relation \eq{eq:wrt}  to $Z_{\CS}$ possible. 
At the same time, it can be shown that a mock theta function and the corresponding false theta function lead to a pair of quantum modular forms that are in fact basically equivalent   (cf. Lemma \ref{lem:asymp}), in a way that precisely leads to the radial limit relation \eq{eqn:rel_CS_orientation}. 
This leads to the natural guess that the quantum invariants $\widehat Z_a(-M_3;\tau)$ for the orientation-reversed three-star plumbed manifold  are given by mock theta functions. This conjecture, proposed in \cite{3d}, will be discussed in \S\ref{mock_conj}.  

In \S\ref{sec:explicitcal}, we will review some recent results supporting the conjecture. The first involves building the relevant orientation-reversed three-manifold via Dehn surgeries on knot complements  \cite{GM}, and the second involves employing the indefinite theta series to extend the definition of $\widehat Z_a(M_3)$ to general plumbed manifolds \cite{future}. To illustrate the various ideas discussed in this proceeding, we will discuss in details the specific example of the Brieskorn sphere $M_3=\Sigma(2,3,7)$. 


\section{False, Mock, and Three-Manifolds} 
\label{sec:false_mock}

In this section we argue that mock modular forms play a role in three-manifold quantum invariants. In \S\ref{subsec:false3star} we introduce the relevant class of quantum invariants and review their relation to false theta functions. In \S\ref{def:falsemockquantum} we review the quantum modular properties of false and mock theta functions and explain their relevance for three-manifold topology. In \S\ref{subsec:mockconj} we discuss a mock conjecture for $\widehat Z_a$ and its motivation and consequences.

\subsection{False Theta Functions and Negative Three-Star Graphs}
\label{subsec:false3star}

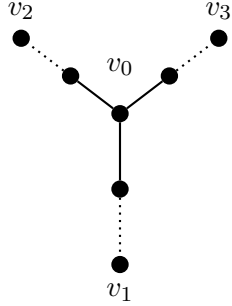
\begin{figure}[H]
	\tikzstyle{point}=[circle,draw=black!100,fill=black!100,thick,
	inner sep=0pt,minimum size=2mm]
	\centering
	\begin{tikzpicture}[thick]
	\node[point] (a0) at ( 0,0)  [label={[label distance=0.25cm] above:{$v_0$}}]{};
	\node[point] (1) at ( 0,-1) {};
	\node[point] (2) at ( -0.65,0.5) {};
	\node[point] (3) at ( 0.65,0.5){};
	\node[point] (a1) at ( 0,-2)[label=below:$v_1$] {};
	\node[point] (a2) at ( -1.3,1)[label=above:$v_2$] {};
	\node[point] (a3) at ( 1.3,1)[label=above:$v_3$] {};
	\draw (1)--(a0)--(2); \draw (a0)--(3); \draw[dotted] (1)--(a1);\draw[dotted] (2)--(a2);\draw[dotted] (3)--(a3);
	\end{tikzpicture}%
	\caption{A three-star graph.}
	\label{fig:3leg}
	\end{figure}

For concreteness, we focus on the simplest type of non-trivial graph: the three-star graphs (see Figure \ref{fig:3leg}). 
These type of graphs correspond via plumbing (cf. 
Figure \ref{fig:plumbing}) to Seifert manifolds with three singular fibers. 
The relation between false theta functions and the WRT invariants for this family of three-manifolds was first pointed out in \cite{LawZag} and later extensively studied in \cite{Hikami1}, \cite{Hikami2}, \cite{Hikami3}. Here we are interested in their quantum invariants 
$\widehat Z_a(M_3)$. It is easy to see \cite{3d} that Definition \ref{def:hatZ} leads to a function well-defined on $\HH$ if and only if $(M^{-1})_{00}<0$, namely when the resulting plumbed three-manifold $M_3$ is weakly negative. 
In order to describe the quantum modular properties of $\widehat{Z}_a(M_3)$, we will need the following definitions. 
\begin{defn}
\label{def:false}
	Let $m\in\mathbb{Z}_{>0}$ and $r \in \mathbb{Z}/2m$. Define false theta functions
\begin{equation}
\widetilde{\theta^1_{m,r}}(\tau):=\sum_{\substack{k \in \mathbb{Z} \\ k \equiv r \text{ } (\text{mod } 2m)}}{\rm sgn}(k)\,q^{\frac{k^2}{4m}}.
\end{equation}
\end{defn}
Note that this is nothing but the usual theta function for one-dimensional lattice $\sqrt{2m}\ZZ$ when the sign factor in the summand is removed. This leads to the nomenclature \cite{RO,FI}. 


It will also be convenient to define, after \cite{plumbing}, the following functions for $m\in\mathbb{Z}_{>0}$ and $j\in \ZZ$
	\begin{equation}\label{def:F}
	F_{j,m}(\tau):=\sum\limits_{k\in \mathbb{Z}}{\rm sgn}\left(k+\frac{1}{2}\right)q^{\left(k+\frac{j}{2m}\right)^2}.
	\end{equation}
 Furthermore,  we have
\begin{equation}\label{rel:Fandfalse}
F_{j,m}(m \tau)=\widetilde{\theta^1_{m,j}}(\tau)+p_{m,j}(\tau) , 
\end{equation}
where $p_{m,j}(\tau)$ is the polynomial in $q$ given by 
\be
p_{m,j}(\tau) = \begin{cases} 
-2\sum\limits_{k=1}^{\lfloor \frac{j}{2m}\rfloor}q^{(j-2mk)^2\over 4m}
&,~{\rm if  }~~ j\geq 2m \\ 
0 & ,~{\rm if  }~~  0\leq j < 2m \\ 
2\sum\limits_{k=0}^{-\lfloor \frac{j}{2m}\rfloor-1}q^{(j+2mk)^2\over 4m} &,~{\rm if  }~~ j < 0.
\end{cases}
\ee
Note that the definition \eq{def:F} can be extended to $m,j\in \QQ^\ast$ since the right hand side only depends on their ratio ${j\over m}$.

In terms of the above building blocks, it can be shown that
given a 3-star weighted graph, the corresponding $\widehat{Z}_a(M_3)$ can be written in terms of $F_{j,p}$ for some $p$ and $j$. See Theorem 4.2 in \cite{plumbing} for the result on a closely related quantity, denoted $Z(q)$ in \cite{plumbing}, and \cite{3d} for numerous examples.

 In particular, in what follows we will further restrict our attention to
 weakly negative plumbed manifolds with three-star plumbing graphs with four nodes. 
 Denote by $M$ the corresponding adjacency matrix, let $A:=-\frac{1}{2}M^{-1}$ and let $v_0$ be the unique vertex with degree three. Moreover, assume that the corresponding adjacency matrix $M$ is unimodular. As a result there is only one quantum invariant $\widehat Z_0(M_3;\tau) :=\widehat Z_a(M_3;\tau) $  with $a=\delta$ mod $2M\ZZ^{|V|}$,  as defined  in 
 \ref{def:hatZ}. Write also
\begin{gather}\label{def:mbc}
\begin{split}
m& = 2 A_{00}\\
    b_0 &= 2\sum_{j=1}^3 A_{j0},  \quad b_i = 4A_{i0}-2\sum_{j=1}^3 A_{j0} \\ 
    c_0 &= A_{12} + A_{23} + A_{31}+ \frac{1}{2}\sum_{j=1}^3 A_{jj}  , \quad c_i = c_0 -
    2 \sum_{\substack{j\in\{1,2,3\}\\ j\neq i}} A_{ij}
\end{split}
\end{gather}
for $i=1,2,3$. Note that $d_i:= -\frac{b_i^2}{4m}+ c_j$ satisfy $d_i=d_j=:d$ for all $i, j\in \{0,1,2,3\}$. 
In the above notation we have the following Proposition. 
\begin{prop}\label{false_prop}
\cite{plumbing}
Consider a weakly negative three-star plumbing graph with four nodes and unimodular adjacency matrix, denote by $M_3$ the corresponding plumbed three-manifold. 
Its unique quantum invariant satisfies
\begin{gather}\begin{split}\label{special_case_false}
(-1)^\pi q^{-c}\widehat{Z}_0(M_{3};\tau)=\sum\limits_{j=0}^3  F_{m- b_j,m}(m\tau)
\end{split}\end{gather}
with $c =d+\frac{3\sigma -\sum_v m(v)}{4}$, where $m$, $b_j$ and $d$ are defined as above and where $\sigma$ and $\pi$ as defined as in Definition \ref{def:hatZ}. 
\end{prop}
Note that, using \eq{rel:Fandfalse} this immediately shows 
\begin{gather}\begin{split}
q^{-c}\widehat{Z}_0(M_{3};\tau)
= \sum\limits_{j=1}^4 \widetilde{\theta_{m,m-b_j}^1}(\tau) + p(\tau) 
\end{split}\end{gather}
where $p(\tau)$ is a polynomial which one can work out explicitly using \eq{rel:Fandfalse}. Often times, one has $-m < b_j \leq m$ for all $j\in\{0,1,2,3\}$ and  $p(\tau)=0$.  
In other words, up to an overall rational power of $q$ and possibly the addition of a polynomial, the quantum invariants $\widehat Z_0$ is given by a false theta function.

\begin{eg}

In this section, we will illustrate the computation of the  quantum invariant and in particular Proposition \ref{false_prop}, with the example of 
 the Brieskorn sphere $M_3=\Sigma(2,3,7)$, which can be described as the intersection between the algebraic surface $\{x^2+y^3+z^7=0\}$ and the five sphere $\{|x|^2+|y|^2+|z|^2=1\}$. It can be obtained as a plumbed manifold with the  plumbing graph shown in Figure \ref{fig:graph237}. 
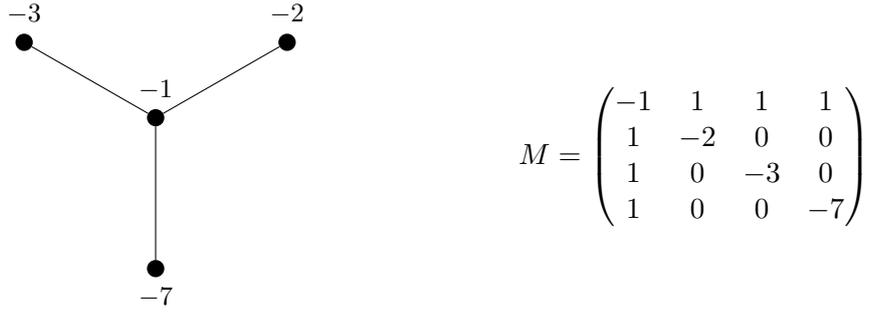
\begin{figure}[h!]
	\tikzstyle{point}=[circle,draw=black!100,fill=black!100,thick,
	inner sep=0pt,minimum size=2mm]
	 \centering
	 \begin{subfigure}[c]{0.4\textwidth}
	 	\centering
	 \begin{tikzpicture}
\node[point] (-1) at ( 0,0)  [label=above:{\small $ -1$}]{};
\node[point] (-7) at ( 0,-2)[label=below:{\small $ -7$}] {};
\node[point] (-3) at (-1.73,1)[label=above:{\small$-3$}] {};
\node[point] (-2) at (1.73,1)[label=above:{\small $-2$}] {};
\draw (-7)--(-1)--(-2); \draw (-1)--(-3);
\end{tikzpicture}%
\end{subfigure}
\begin{subfigure}[c]{0.4\textwidth}
	\begin{equation*}
	M=
	\begin{pmatrix}
	-1 & 1& 1& 1\\
1&	-2 & 0 & 0  \\
1&	0 & -3 & 0  \\
 1&	0 & 0 & -7 \\
	\end{pmatrix}
	\end{equation*}
\end{subfigure}
	\caption{Plumbing graph and adjacency matrix for $\Sigma(2,3,7)$}
\label{fig:graph237}
\end{figure}%
Note that $M$ is indeed unimodular, consistent with the fact that Brieskorn sphere is a integral homology sphere with trivial $H_1(M_3;\ZZ)\cong \ZZ^{4}/M\ZZ^{4}$ and there is hence just one quantum invariant $\widehat{Z}_0(M_3;\tau)$. 

Plugging the adjacency matrix in \eq{def:mbc} one obtains 
\be\label{mbc237}
m=42~,~~(b_j,4c_j)= (1,1),(-13,5),(-29,21),(41,41)~{\rm for}~j=0,1,2,3. 
\ee
Using Definition \ref{def:hatZ} (or \eq{special_case_false}), one obtains that 
\begin{gather}
\label{eqn:237example}\begin{split}
\widehat{Z}_0(\Sigma(2,3,7);\tau) 
&= q^{\frac{83}{168}} \sum_{\substack{k\geq 0 \\ k^2 \equiv 1 ~~(42)}} \left(\frac{k}{21}\right) q^{\frac{k^2}{168}} \\
&= q^{\frac{83}{168}}\left(\widetilde{\theta_{42,1}^1} -\widetilde{\theta_{42,13}^1} -\widetilde{\theta_{42,29}^1} +\widetilde{\theta_{42,41}^1}    \right) (\tau)\\
&= q^{\frac{1}{2}} \left(1-q-q^5+q^{10}-q^{11}+\dots \right) 
\end{split}\end{gather}
which is indeed a false theta function of weight 1/2. 
The fact that $Z_{\CS}(\Sigma(2,3,7))$ is given by the above function by taking the limit 
\be\label{radial_limit_237}
Z_{\CS}(\Sigma(2,3,7);k) ={1\over i\sqrt{2k}} \lim_{t\to 0^+} (\widehat{Z}_0(\Sigma(2,3,7);\frac{1}{k}+t))
\ee
was first established by \cite{LawZag}. 
\end{eg}

\subsection{False, Mock, and Quantum}
\label{def:falsemockquantum}

As we have seen, a pre-requisite for a $q$-series to play the role of the quantum invariants $\widehat Z_a(M_3)$ is to have a specific behaviour when taking the radial limit, so that it gives the WRT invariants via \eq{eq:wrt}. 
This is  demonstrated in the $\Sigma(2,3,7)$ example in \eq{radial_limit_237}. This leads us to the concept of {\em quantum modular forms} (QMF), first introduced by D. Zagier  \cite{qmf}.
Roughly speaking, a quantum modular form is a function defined on $\QQ$ with a certain modular-like property: the deviation from modularity, measured by a modular difference function denoted by $p_\gamma$, has nice analytic properties that are not a priori manifest or expected. In this proceeding we work with a specific version of the definition that is often referred to as defining  {\em strong quantum modular forms}. We refer to \S7.3 of \cite{3d} for details about  modular forms in the current context.  

To state the definition, let us recall the usual definition 
of the {slash operator}  acting on the space of holomorphic functions on $\mathbb{H}$ for weight $w$ and multiplier $\chi$ on $\Gamma$, which we take to be a subgroup of the modular group $\SL_2(\ZZ)$: 
\be\label{def:slash}
f(\tau)\lvert_{w,\chi} \gamma := f\left(\frac{a\t+b} {c\t+d}\right)  \chi(\gamma) (c\t+d)^{-w}\,, \qquad \gamma=\left( \begin{smallmatrix} a&b\\c&d\end{smallmatrix}\right)\in \Gamma. 
\ee


\begin{defn} 
   \label{def:qu}\cite{qmf} 
   Consider a function $Q: \QQ \rightarrow \CC$. It is called 
a \emph{strong quantum modular form} of weight $w$ and multiplier $\chi$ for $\Gamma$ if for every $\gamma\in \Gamma$ the modular difference function $p_\g (x): \QQ \backslash \{\g^{-1}(\infty)\} \rightarrow \CC$, defined by 
\be
\label{def:period}
p_\g (x):=Q(x) - Q\lvert_{w,\chi} \gamma(x)
\ee
is a real-analytic function of $\RR$ minus finitely many points. 
   \end{defn}

The false theta functions we encountered in \S\ref{subsec:false3star} are examples of quantum modular forms. The simplest way to see this is to note that false theta functions defined in \eq{def:false} are examples of Eichler integrals. Given a cusp form $g=\sum_{n>0} a_g(n) q^n$ of weight $w \in \frac{1}{2}\ZZ$, its Eichler integral is defined as
\be\label{def:eichler}
\widetilde{g}(\tau) : =\sum_{n>0} n^{1-w} a_g(n) q^n .
\ee
It is easy to see that the false theta function $\widetilde{\theta_{m,r}^1}$ is the Eichler integral of the weight 3/2 unary theta function 
\be\label{unary theta}
\theta^1_{m,r}(\tau)  : ={1\over \sqrt{4m}}\sum_{\substack{\ell\in\ZZ\\\ell=r\xmod 2m}}\ell \, q^{\ell^2/4m}, 
\ee
as the notation suggests. 

Note that this is equal to the following integral\footnote{We choose the branch to be the principal branch $-\pi < {\rm arg}x \leq \pi$. }
\be\label{def:Eichler_int1}
\widetilde{g}(\tau) = C \int_{\tau}^{i\infty} g(z' ) (z'-\t)^{w-2} dz',
\ee
where $C=\frac{(2\pi i)^{w-1}}{\Gamma(w-1)}$.
Letting $Q(x):=\lim_{t\to 0^+} \widetilde{g}(x+it)$ in Definition \ref{def:qu}, one immdiately sees that the modular difference function $p_\g (x)$ admits an expression as a period integral 
\be\label{periodEich}
p_\gamma(x) = C  \int_{\gamma^{-1}(i\infty)}^{i\infty} g(z' ) (z'-x)^{w-2} dz'
\ee
and  is hence equipped with the desired analytic properties. 

An analogous argument demonstrates that mock modular forms also lead to quantum modular forms.
To see that,  we first recall the definition of mock modular forms, adapted to the classes of functions that are relevant for our context. In particular, we assume that the shadow $g$ is a cusp form with real Fourier coefficients, namely $\overline{g(-\bar \tau)}=g(\tau)$.

\begin{defn}\label{def:mock}
We say that a holomorphic function $f$ on $\HH$ is a {\em mock modular form} of weight $k$ and multiplier $\chi$ on $\Gamma$ if there exists a weight $2-k$ cusp form $g$ on $\Gamma$
such that the non-holomorphic {\em completion} of $f$, defined as
$$\hat f (\t) = f(\tau) - g^\ast (\tau) $$
satisfies $\hat f = \hat f\lvert_{k,\chi} \gamma$ for every $\gamma\in \Gamma$.
In the above, $g^\ast$ denotes the non-holomorphic Eichler integral
\be
g^\ast(\t) 
:=
 { C\int_{-\bar \t}^{i\infty} (\t' + \t)^{-k} {g(\t')} \, d\t' },
\ee
defined for $\t\in \HH$.
\end{defn}
Note that there is no canonical normalization for the shadow and  we  choose ours to simplify the comparison between mock modular forms and Eichler integrals \eq{def:eichler}. \underline{Assuming} that the limit  $\lim_{t\to 0^+} f(x+it)$ exists for a given $x\in \QQ$, let $Q(x):= \lim_{t\to 0^+} f(x+it)$,  
and consider $\gamma\in \Gamma$ in the notation of Definition \ref{def:qu}. The $\gamma$-invariance of the completion $\hat f$ leads to an expression for the corresponding modular difference function $p_\gamma(x)$  given by the modular difference function associated to $g^\ast(\tau)$. 
Through the related function (given by $\tau\mapsto -\tau$) 
\be
\tilde g^\ast(\tau)=
 { C\int_{\bar \t}^{i\infty} (\t' - \t)^{-k} {g(\t')} \, d\t' }
\ee
and the fact that the modular difference function associated with $\tilde g^\ast$ is again
a period integral completely analogous to \eq{periodEich}, it follows that 
mock modular forms indeed lead to quantum modular forms. 
Moreover, the quantum modular forms arising from a mock modular form and the Eichler integral of its shadow are clearly closely related. 


However, just as Ramanujan already pointed out in his original work \cite{Ram88}, mock modular forms inevitably encounter infinities when approaching certain rational numbers from within the upper-half plane. Nonetheless, there exists a finite collection  of weakly holomorphic modular forms that can be used to ``cut out" these infinities and render a well-defined radial limit. More explicitly, we have the following theorem. 

\begin{thm}\label{thm_mock_Ram} \cite{MR3065809, CLR}
Let $f$ be a mock modular form of weight $k$ and multiplier system $\chi$ for $\Gamma_0(N)$ with non-vanishing shadow $g$, and let $\{x_1,\dots, x_t\}\subset \QQ\cup \{i\infty\}$ be a set of representatives of $\Gamma_0(N)$-inequivalent cusps, 
then\vspace{-5pt}
\begin{enumerate}
\item the function $f(\t)$ has exponential singularities at infinitely many rational numbers, \vspace{-7pt}
\item for every weakly holomorphic modular form $G$ of weight $k$ and multiplier system $\chi$ for $\Gamma_0(N)$, $f-G$ has exponential singularities at infinitely many rational numbers, \vspace{-7pt}
\item there exits a collection $\{G_j\}_{j=1}^{t}$ of weakly holomorphic modular forms with the following property. Given any cusp $x$, let $x_j$ be the cusp representative that is $\Gamma_0(N)$-equivalent to $x$ and write $G_x= G_j$. Then 
$f-G_x$ is bounded towards $x$.
\end{enumerate}
\end{thm}

Moreover, following the arguments sketched above, the mock modular form and the Eichler integral of its shadow leads to a pair of closely related strong quantum modular forms. 
\begin{lem}\label{lem:asymp}
\cite{3d}
With the notation of Theorem \ref{thm_mock_Ram}, let $g$ be the shadow of $f$, 
the asymptotic expansions of the Eichler integral $\tilde g$ and the mock modular fomr $f-G_x$ near $x$ take the form
\be
 (f-G_x)(-x+it) \sim  \sum_{n\geq 0} \alpha_{x}(n) (-t)^n ~{\rm and}~\tilde g(x+it)\sim \sum_{n\geq 0}
 \alpha_{x}(n) t^n.
 \ee
\end{lem}

In particular, when the shadow $g$ is a weight 3/2 unary theta function, the mock modular forms are (up to an overall rational power of $q$) called mock theta functions in the terminology of 
\cite{zagier_mock}, and the Eichler integral are the false theta functions encountered in \S\ref{subsec:false3star}. The false-mock pair satisfies
\be
\label{quantum_mock}
\lim_{t\to 0^+} (f-G_x)(x+it) = \lim_{t\to 0^+} \tilde g(-x+it), 
\ee
reminiscent of the relation \eq{eqn:rel_CS_orientation} between $Z_{\CS}(M_3)$
 and  $Z_{\CS}(-M_3)$ when taking $x=\frac{1}{k}$. 
Focusing on the cusp $x=0$, we can see that the false and mock forms have the ``same'' asymptotic series, approaching from the upper- and lower-half plane, in the sense that the asymptotic expansions in the limit $t\to 0^+$ satisfy
\be \label{expansionzero}
 (f-G_0)(it) \sim  \sum_{n\geq 0} \alpha_{0}(n) (-t)^n ~{\rm and}~\tilde g(it)\sim \sum_{n\geq 0}
 \alpha_{0}(n) t^n .
\ee
On the three-manifold side, the cusp $x=0$ is relevant for the perturbative invariants (the so-called Ohtsuki series), capturing the expansion around the semi-classical $k\to \infty$ limit. 

To end this subsection, we provide an explicit example of such a false-mock pair. 
\begin{eg}
Consider the order seven mock theta function $F_0(q)$ by Ramanujan\cite{Ram00}. It is, up to an overall power of $q^{-\frac{1}{168}}$, a mock modular form of weight 1/2
\be
f(\tau)= {q^{-\frac{1}{168}}}\,F_0(q)  = {q^{-\frac{1}{168}}} \left(1+q+q^3+q^4+ q^5+2q^7+O(q^8)\right), 
\ee
whose shadow is given by the unary theta function
\be\label{237unarytheta}
g(\tau) = \left(\theta^1_{42,1}-\theta^1_{42,13}-\theta^1_{42,29}+\theta^1_{42,41}\right)(\tau). 
\ee
Compared to \eq{eqn:237example}, we see that the Eichler integral is (up to a factor $q^{\frac{83}{168}}$) precisely the 
quantum invariant of the Brieskorn sphere $\Sigma(2,3,7)$:
\be
\widehat Z_0(\Sigma(2,3,7);\tau) = q^{\frac{83}{168}} \tilde g(\tau). 
\ee
\end{eg}

\subsection{A Mock Conjecture}
\label{subsec:mockconj}

In \cite{3d}, the following relation between mock modular forms and three-manifold quantum invariants is proposed\footnote{Note that when $-M_3$ is not a weakly negative plumbed manifold, the mathematical definition \ref{def:hatZ} does not apply and this conjecture can be seen as rather a definition. However, recall that a physical definition of $\widehat Z_a(M_3)$ does exist for all closed three-manifolds \cite{BPS}. As a result,  independent computations can  in principle be carried out for $-M_3$, as we will  demonstrate in  \S\ref{subsec:negknots} for certain classes of $-M_3$. With this in mind we regard \eq{conj:mM3} as a conjecture.   }. 
\begin{conj}\label{mock_conj}
Let $M_3$ be a three-manifold whose quantum invariants take the form 
\be
\widehat Z_a(M_3;\tau) = q^{c} \left( \tilde g(\tau) + p(\tau)\right) 
\ee
where $c\in \QQ$, $\tilde g(\tau)$ is the Eichler integral of a theta function $g(\tau)$ of weight $w=\frac{3}{2}$ and $p(\tau)$ is a polynomial in $q$, then
\be\label{conj:mM3}
\widehat Z_a(-M_3;\tau) = q^{-c} \left( f(\tau) + p(-\tau)\right), 
\ee
where $f(\tau)$ is a weight $1/2$ mock modular form whose shadow is given by $g(\tau)$. 
\end{conj}

The relevance of the above conjecture can be seen in Proposition \ref{false_prop}, which guarantees the existence of $M_3$ satisfying the condition of the conjecture. More generally, we also expect mixed weight and higher-depth mock modular forms to play a role in three-manifolds quantum invariants. See \cite{3d} and \cite{3dlogvoa}.
In what follows we briefly describe the 
three general motivations for the above conjecture, first discussed in \cite{3d}. In \S\ref{sec:explicitcal} we will present explicit calculations which render results predicted by Conjecture \ref{mock_conj}, and hence constitute further evidence for it. 
\begin{itemize}
\item As mentioned in the previous subsection, the asymptotic values \eq{quantum_mock} and expansions \eq{expansionzero} of a false-mock pair are analogous to the relation  \eq{eqn:rel_CS_orientation} among the WRT invariants of a pair of three-manifolds related by a flip in orientation.  

\item Some false theta functions have known expressions as $q$-hypergeometric series, which converge not only inside but also outside the unit circle (when considered as a function of $q$). In some cases  the expression on the other side is given by a mock theta function. See \S7.4 of \cite{3d} for details. 

\item When a weight $1/2$ mock modular form can be expressed as a so-called Rademacher sum, one can prove in general that the same Rademacher sum, now performed in the lower rather than the upper half-plane, yields precisely the corresponding Eichler integral. In other words, the Rademacher sum yields a function defined on both $\HH$ and $\HH^-$, where they coincide with the mock resp. false theta function in question.

\end{itemize}

We refer to \S7.4 of \cite{3d} for a detailed discussion of the third point above. 
To illustrate the second point, let us  consider an example that is again relevant for the Brieskorn sphere $M_3=\Sigma(2,3,7)$. \begin{eg}
Let us  define a function $\psi: \HH\cup \HH^-$ in terms of the $q$-hypergeometric series: 
\be
\psi(\tau) := \sum_{n\geq 0}\frac{ q^{n^2}}{(q^{n+1};q)_{n}}. 
\ee
 Note that the $q$-hypergeometric series converges both for $|q|<1$ and $|q|>1$. 
It can be shown that \cite{Hikami2}
\be
\psi(\tau) = \begin{cases} q^{-\frac{1}{168}}\tilde g(\tau) &, ~~\tau \in \HH \\ 
 F_0(q^{-1}) &, ~~\tau \in \HH^- 
\end{cases}. 
\ee
See also \cite{BFR} for a more general discussion.
As a result, since $ \widehat Z_0(\Sigma(2,3,7),\tau)=q^{\frac{1}{2}} \psi(\tau)$ for $\tau\in \HH$, we can try to extend the definition of LHS to $\HH^-$ using the RHS. 
It is hence natural to guess that (cf. \eq{eqn:reflect})
\be \label{guess_237}
\widehat Z_0(-\Sigma(2,3,7),\tau) ~{\rm ``="}~ \widehat Z_0(\Sigma(2,3,7),-\tau) {\rm ``="}  q^{-\frac{1}{2}}  F_0(q). 
\ee
\end{eg}

We now end this section with a discussion on certain important open questions. 
First, note that Conjecture \ref{mock_conj} does not specify, given a shadow, {\em which} mock modular form $f$ should be the correct quantum invariant for the orientation-reversed manifold $M_3$. Recall that two mock modular forms differing by a (weakly holomorphic) modular form have the same shadow.  
This question is of crucial importance since, as proposed in \cite{BPS}, 
the Fourier coefficients of the quantum invariants $\widehat Z_a$ are (up to a possible factor of 2) integers which have the physical interpretation of counting supersymmetric quantum states in the underlying quantum physical theory. This said, we do expect the leading term of $\widehat Z_a$ in the $\tau\to i\infty$ expansion to obey the naive $q\leftrightarrow q^{-1}$ relation and this puts meaningful constraints on the mock modular forms. 
Second, as we have seen in \S\ref{def:falsemockquantum}, mock and false theta functions relate to the WRT invariants in a slightly different way. While the radial limit of false theta functions are well-defined, for many cusps $x$ one has to subtract the singular terms (by subtracting a modular form $G_x$ which cuts out the singularity for instance) of the mock form in order to have a well-defined limit when approaching $x$ from within the upper-half plane (cf. Lemma \ref{lem:asymp}). The asymmetry might not be so surprising from the physical point; the $M5$-brane theory is known to be a chiral theory. It would be extremely interesting to understand the physical or topological interpretation of the singular terms when taking radial limit of mock theta functions.

\section{Explicit Calculations}
\label{sec:explicitcal}

In this section we summarise recent developments which make it possible to define and to compute the quantum invariants  $\widehat Z_a(-M_3)$ for certain three-manifolds $-M_3$ that are relevant for the mock conjecture discussed in  \S\ref{subsec:mockconj}. 
 We illustrate these methods with explicit computations for the Brieskorn sphere $M_3=\Sigma(2,3,7)$.
 
\subsection{Quantum Invariants via Knots}
\label{subsec:negknots}

In this subsection, we review a (conjectural) way, introduced in \cite{GM}, to compute the quantum invariants $\widehat Z_a$ for some of the three-manifolds  that are relevant for the mock conjecture \ref{mock_conj}, 
by constructing them via Dehn surgeries of knot complements.



Consider a knot $K$. Let $Y(K)$ be the knot complement of $K$ in an integral homology sphere $\hat Y$. A closed manifold $Y_{p/r}(K)$ can be obtained by $Y(K)$ via Dehn surgery with coefficient $p/r\in \QQ^\ast$. Roughly speaking, $p/r$ specifies the diffeomorphism of $\partial Y(K)$, dictating the way a solid torus is glued along  $\partial Y(K)$ to obtain  $Y_{p/r}(K)$.  

Now consider the special case when $\hat Y=S^3$. 
Given this choice, one associates to a knot $K$ a  two variable series 
\begin{equation}
F_K(x,q)\in 2^{-c}q^{\Delta}\mathbb{Z}[x^{1/2,},x^{-1/2}][q^{-1},q]]
\end{equation}
where $c\in \mathbb{Z}_{+}$ and $\Delta\in \mathbb{Q}$.
For instance, for $K$ a positive torus knot, an explicit expression for $F_K(x,q)$  has been given in \cite{GM}. 
Define a ``Laplace transform" $\mathcal{L}^{(a)}_{p/r}$,  given by (see also \cite{GMP})
\begin{equation}
\mathcal{L}^{(a)}_{p/r}: x^uq^v\to \begin{cases}
\begin{array}{ll}
q^{-u^2r/p}\cdot q^v & \text{if } ru-a\in p\mathbb{Z}\,, \\
0 & \text{otherwise}\,.
\end{array}
\end{cases}
\end{equation}
It has been shown for positive torus knots $K$ (Theorem 1.2 of \cite{GM}) and conjectured for general knots (Conjecture 1.7 of \cite{GM}) that, for values of $p/r$ such that the right hand side is well defined and for some $d \in \mathbb{Q}$ and $\varepsilon\in\{\pm1\}$, one has
\begin{equation}
\label{resurgenceinvariants}
\widehat{Z}_a(\tau, S^3_{p/r}(K))=\varepsilon q^d\cdot \mathcal{L}^{(a)}_{p/r}[(x^{\frac{1}{2r}}-x^{-\frac{1}{2r}})F_K(x,q)]
\end{equation}
 where we canonically identify the ${\rm Spin}^c$-structure of $S^3_{p/r}(K)$ with 
\be
a\in \ZZ+\frac{r+1}{2} \xmod p\ZZ. 
\ee

Now it remains to compute $F_K$ for general knots. 
It is convenient to define a rescaled version of $F_K(q,x)$: 
\begin{equation}
f_K(x,q):=\frac{F_K(x,q)}{x^{\frac{1}{2}}-x^{-\frac{1}{2}}}.
\end{equation}
Based on physical expectations, a relation between the Borel resummation of the colored Jones polynomial $J_n(e^{\hbar})$ and $f_K(x,q)$, where $q=e^\hbar$ and $x=e^{n\hbar}$, is conjectured in (Conjecture 1.5 of) \cite{GM}. Drawing inspiration from the analogous conjectures for the Chern-Simons partiton function on the knot complement \cite{Guk} or for the colored Jones polynomials \cite{Gar}, the following was proposed in \cite{GM}
\begin{conj} For any knot $K \subset S^3$ the quantum polynomial $\widehat{A}$ of $K$ annihilates the series $f_K(x,q)$
\begin{equation}\label{recursion}
\widehat{A}f_K(x,q)=0
\end{equation}
and
\begin{equation}
\label{boundary}
\lim_{q \to 1}f_K(x,q)=s.e.\frac{1}{\Delta_K(x)}
\end{equation}
where the symmetric expansion s.e. denotes the average of the expansions of the given rational function as $x \to 0$ (as a Laurent power series in x) and as $x \to \infty$ (as a Laurent power series in
$x^{-1}$).
\end{conj}
Note that  \eqref{recursion} sets up a recursion relation for the coefficients $f_m(q)$ in $f_K(q,x)=\sum\limits_{m}f_m(q)x^m$, while the relation \eq{boundary} to the Alexander polynomial $\Delta_K(x)$ provides a boundary condition for the recursion equation.
This is often sufficient to determine $F_K$ to any desired order. 

\begin{eg}
For the figure-eight knot $K=\pmb{4_1}$, the above-mentioned procedure leads to the leading order expansion\cite{GM}
\be\label{eqn:figure81}
F_{\pmb{4_1}}(x,q) ={1\over2}\left(\Xi(x,q)-\Xi(x^{-1},q)\right) 
\ee
where
\be\label{eqn:figure82}
\Xi(x,q) =x^{1/2}+2\,x^{3/2}+(q^{-1}+3+q)x^{5/2} + 
(2q^{-2}+2q^{-1}+5+2q^{1}+2q^{2}) +\dots . 
\ee

The orientation-flipped Brieskorn sphere $-\Sigma(2,3,7)$ can be constructed through surgery on the complement in $S^3$ of the figure-eight knot \pmb{$4_1$}, namely $-\Sigma(2,3,7)=S^3_{-1}(\pmb{4_1})$.
Exploiting the conjecture \eq{resurgenceinvariants} and plugging in \eq{eqn:figure81}-\eq{eqn:figure82}, we obtain the result:
\begin{equation}\label{dehn237}
\widehat{Z}_0(-\Sigma(2,3,7))=-q^{-\frac{1}{2}}(1 + q + q^3 + q^4 + q^5 + 2q^7 + q^8 + 2q^9 + q^{10} + 2q^{11} + \ldots). 
\end{equation} 
Note that the above leading terms in the $q$-expansion coincide (up to a sign) with the guess \eq{guess_237} based on quantum modular properties and on $q$-hypergeometric identities. However, the procedure outlined in this subsection does not immediately lead to a way to prove the modularity of \eq{dehn237}. We will see yet another way to compute $\widehat{Z}_0(-\Sigma(2,3,7))$  in the following subsection.

\end{eg}

\subsection{Relation to indefinite theta functions}
\label{subsec:indefinite}

As mentioned earlier, one immediate problem with the 
proposal $\widehat Z_a(-M_3;\tau) {\text{ ``=''}}\, \widehat Z_a(M_3;-\tau)$ \eq{eqn:reflect}
is the fact that in Definition \ref{def:hatZ} one has 
$\tau \leftrightarrow -\tau \Leftrightarrow M \leftrightarrow -M$, and after this flipping of  signature one no longer obtains a  theta function $\Theta^{-M}_a$ \eq{def:theta} (and an integral \eq{def:hatZ}) that makes sense on the upper-half plane.  


While it seems to be the end of the road as far as  Definition \ref{def:hatZ} is concerned,  a natural possibility is to replace the  naive theta series   with a regularised theta function. Indeed, building on previous work by 
Vign\'eras \cite{Vigneras}, Zwegers \cite{zwegers} has devised a way to define a regularisation for theta functions of signature $(1,n)$ which retains its holomorphicity, and moreover established the relation to mock theta functions. The regularisation of  general indefinite theta functions and the relation to higher-depth mixed mock modular forms has recently been developed in \cite{ABMP}, \cite{Naz}, \cite{Wes}, \cite{FK}, \cite{ZZ}. In \cite{future}, we apply
 these results to define and to compute quantum invariants for plumbed three-manifolds 
 that are not weakly negative.  

For the sake of concreteness and in order to establish a direct relation to Ramanujan's mock theta function, we focus on the class of three-manifolds discussed 
in Proposition \ref{false_prop}. 
In the notation of Proposition \ref{false_prop} and of Figure \ref{fig:3leg}, after performing the integration over $w_{v_i}$ for $i\in \{1,2,3\}$ and write the $w_{v_0}=w$, we obtain 
\be
\widehat Z_0(M_3;\tau) = (-1)^{\pi}q^{\frac{3\sigma -\sum_v m(v)}{4}}\, \text{vp} \oint_{|w|=1}\frac{dw}{2\pi i w (w-w^{-1})} h(\tau, z)
\ee
where 
\be
 h(\tau, z) = \sum_{j=0}^3 \sum_{\varepsilon \in \{\pm 1\}} \varepsilon \sum_{k\in 1+2\ZZ} q^{\frac{m}{4}k^2 -\frac{\varepsilon b_j}{2}k+c_j} w^k . 
\ee
Note that naively taking $\tau \mapsto -\tau$ in $ h(\tau, z)$ gives
\begin{gather}\begin{split}
&\frac{q^{\frac{1}{24}}}{\eta(\tau)} \sum_{j=0}^3 \sum_{\varepsilon \in \{\pm 1\}} \varepsilon \sum_{k\in 1+2\ZZ} \sum_{n\in \ZZ} (-1)^n q^{-\frac{m}{4}k^2 +\frac{\varepsilon b_j}{2}k-c_j + \frac{3n^2-n}{2}} w^k \\
&= \frac{q^{-d}\, e^{\frac{\pi i}{6}} }{\eta(\tau)}\sum_{j=0}^3 \sum_{\varepsilon \in \{\pm 1\}} \varepsilon  w^\frac{\varepsilon b_j}{m} \sum_{{\bf v} \in \Lambda_{j,\varepsilon}}q^{\frac{({\bf v},{\bf v})}{2}}e^{2\pi i (z ~\frac{1}{2}).{\bf v}}
\end{split}
\end{gather}
where we have inserted $1={\eta(\tau)\over \eta(\tau)} = \frac{{q^{\frac{1}{24}} \sum_{n\in \ZZ} (-1)^n q^{ \frac{3n^2-n}{2}}}}{\eta(\tau)}$, and the bilinear form in the second line is given by 
\be
({\bf v}',{\bf v}):= {{\bf v}^T K {\bf v}}, ~~ K:= \bem -\frac{m}{2} & 0 \\ 0&3\eem,
\ee
and the set of summation is given by 
\be
\Lambda_{j,\varepsilon}= \left\{
{{\bf v} = (\begin{smallmatrix} v_1\\ v_2 \end{smallmatrix}) \large\lvert v_1 \in 2\ZZ +1 -\frac{\varepsilon b_j}{m} ,v_2 \in \ZZ-\frac{1}{6}}
\right\}.
\ee

In other words, the key ingredient $h(\tau,z)$ of the integrand becomes, after taking $\tau \mapsto -\tau$ and multiplying by $\eta(\tau)$,  a sum of theta functions of signature $(1,1)$ that we would like to make sense of. 
As a result, we propose the following definition for three-manifolds $M_3$ that satisfy the conditions specified in Proposition \ref{false_prop} 
({\it i.e.}, when $M_3$ is a weakly negative plumbed manifold which can be obtained from  a  three-star plumbing graph with four nodes and unimodular adjacency matrix):
\be\label{def:negativehatZ}
\widehat Z_0(-M_3;\tau) := \frac{ (-1)^{\pi}q^{-\frac{3\sigma -\sum_v m(v)}{4}}}{\eta(\tau)}
\text{vp} \oint_{|w|=1}\frac{dw}{2\pi i w (w-w^{-1})} \vartheta^M_a(\tau, z)
\ee
where 
\begin{gather}\begin{split}
\vartheta^M_a(\tau, z) &:= {q^{-d}\, e^{\frac{\pi i}{6}} }\sum_{j=0}^3 \sum_{\varepsilon \in \{\pm 1\}} \varepsilon  w^\frac{\varepsilon b_j}{m} \sum_{{\bf v} \in \Lambda_{j,\varepsilon}}\rho({\bf v})\,q^{\frac{({\bf v},{\bf v})}{2}}e^{2\pi i (z ~\frac{1}{2}).{\bf v}} 
\end{split}\end{gather}
for an appropriately chosen ``regularisation factor" $\rho({\bf v})$ which will be discussed in details in \cite{future} and 
will be described explicitly in the example below.

\begin{eg} We will again take the example of $M_3 =\Sigma(2,3,7)$,  with the plumbing graph and the adjacency matrix given in Figure \ref{fig:graph237}. The relevant parameters $m$, $b_j$, and $c_j$ are given in \eq{mbc237}.
Adapting \cite{zwegers} to preserve the symmetry\footnote{In this case, as noted in \cite{plumbing}, the principal value contour integral renders the same result as integrating over a contour lying inside (or outside) the unit disk.} $\vartheta^M_a(\tau, z) = - \vartheta^M_a(\tau, -z)$, we choose the regularising factor
\be\label{eqn:reg_237}
\rho({\bf v})=\rho^{c,c'}({\bf v}): = \frac{1}{2} \left({\rm sgn}(\bar {\bf v}, c) -{\rm sgn}(\bar {\bf v}, c')\right) 
\ee
where $\bar {\bf v} = (\begin{smallmatrix} |v_1|\\ v_2 \end{smallmatrix}) $ for ${\bf v} = (\begin{smallmatrix} v_1\\ v_2 \end{smallmatrix})$ and 
the two timelike vectors are chosen to be 
\be\label{eqn:choiceofc} c=(1,0), ~~ c'= (8,21).\ee See \cite{future} for a discussion on the choices of the timelike vectors $c$ and $c'$. 
Putting things together, we have the following result (see \cite{future} for a proof).  \end{eg}
\begin{prop}
When the regularisation factor $\rho({\bf v})$ is given as in \eq{eqn:reg_237}, \eq{eqn:choiceofc}, the definition \eq{def:negativehatZ} leads to 
\begin{gather}\begin{split}
q^{\frac{1}{2}} \widehat Z_0(-\Sigma(2,3,7);\tau) &= F_0(q)\\
&=1+q+q^3+q^4+ q^5+2q^7+O(q^8)
\end{split}\end{gather}
where $F_0(q)$ is the order 7 mock theta function of Ramanujan. 
\end{prop}

Note that this result, given by the order 7 mock theta function, was precisely what was expected in \cite{3d} \eq{guess_237}. 
Moreover, at least the leading orders of $q$-expansion also, up to a sign, coincides with the result \eq{dehn237} which was obtained via a logically totally independent computation. These results constitute supporting evidence for Conjecture \ref{mock_conj}.

\section*{Acknowledgements}
We would like to thank Sungbong Chun, Sergei Gukov and Sarah Harrison for extremely helpful conversations.
The work of F.F. is supported in part by the MIUR-SIR grant RBSI1471GJ ``Quantum Field Theories at Strong Coupling: Exact Computations and Applications". The work of M.C. and G.S. is supported by the NWO vidi grant (number
016.Vidi.189.182). The work of M.C. has also received support from ERC starting grant H2020 \#640159.

\end{document}